\documentclass[10pt,letterpaper,dvips]{article}

 \usepackage{amssymb,latexsym,amsmath,amsthm}
\newtheorem{thm*}{Theorem}

\newtheorem{thm}{Theorem}[section]

\newtheorem{dfn}{Definition}[section]

\newtheorem{lemma}{Lemma}[section]

\newtheorem{cor}{Corollary}[section]

\newtheorem{prop}{Proposition}[section]

\usepackage{fullpage}

\begin{document}

\def\d{ \partial }
\def\Na{{\mathbb{N}}}

\def\Z{{\mathbb{Z}}}

\def\IR{{\mathbb{R}}}

\def\L{ {\mathcal{L}}}

\newcommand{\E}[0]{ \varepsilon}

\newcommand{\s}[0]{ \mathcal{S}}

\newcommand{\AO}[1]{\| #1 \| }

\newcommand{\BO}[2]{ \left( #1 , #2 \right) }

\newcommand{\CO}[2]{ \left\langle #1 , #2 \right\rangle}

\newcommand{\co}[1]{ #1^{\prime}}

\newcommand{\p}[0]{ p^{\prime}}

\newcommand{\m}[1]{   \mathcal{ #1 }}

\newcommand{ \A}[1]{ \left\| #1 \right\|_H }

\newcommand{\B}[2]{ \left( #1 , #2 \right)_H }

\newcommand{\C}[2]{ \left\langle #1 , #2 \right\rangle_{  H^* , H } }

 \newcommand{\HON}[1]{ \| #1 \|_{ H^1} }

\newcommand{ \Om }{ \Omega}

\newcommand{ \pOm}{\partial \Omega}

\newcommand{\D}{ \mathcal{D} \left( \Omega \right)}

\newcommand{\Ov}{ \overline{ \Omega}}

\newcommand{\DP}{ \mathcal{D}^{\prime} \left( \Omega \right)  }

\newcommand{\DPP}[2]{   \left\langle #1 , #2 \right\rangle_{  \mathcal{D}^{\prime}, \mathcal{D} }}

\newcommand{\PHH}[2]{    \left\langle #1 , #2 \right\rangle_{    \left(H^1 \right)^*  ,  H^1   }    }

\newcommand{\PHO}[2]{  \left\langle #1 , #2 \right\rangle_{  H^{-1}  , H_0^1  }}

 \newcommand{\HO}{ H^1 \left( \Omega \right)}

\newcommand{\HOO}{ H_0^1 \left( \Omega \right) }

\newcommand{\CC}{C_c^\infty\left(\Omega \right) }

\newcommand{\N}[1]{ \left\| #1\right\|_{ H_0^1  }  }

\newcommand{\IN}[2]{ \left(#1,#2\right)_{  H_0^1} }

\newcommand{\INI}[2]{ \left( #1 ,#2 \right)_ { H^1}}

\newcommand{\HH}{   H^1 \left( \Omega \right)^* }

\newcommand{\HL}{ H^{-1} \left( \Omega \right) }

\newcommand{\HS}[1]{ \| #1 \|_{H^*}}

\newcommand{\HSI}[2]{ \left( #1 , #2 \right)_{ H^*}}

\newcommand{\WO}{ W_0^{1,p}}
\newcommand{\w}[1]{ \| #1 \|_{W_0^{1,p}}}

\newcommand{\ww}{(W_0^{1,p})^*}

\newcommand{\labeld}[1]{ \mbox{ \qquad (#1)}  \qquad   \label{ #1}  }

\newcommand{\bi}{\Delta^2}
\newcommand{\la}{\lambda}
\newcommand{\f}{\frac}
\newcommand{\wt}{\widetilde}
\newcommand{\dis}{\displaystyle}
\newcommand{\om}{\Omega}

\title{The critical dimension for a $4^{\rm th}$ order problem with singular nonlinearity}

\author{{Craig Cowan\footnote{Department of Mathematics, University of British Columbia, Vancouver, B.C. Canada V6T 1Z2. E-mail: cowan@math.ubc.ca. }
\qquad
Pierpaolo Esposito\footnote{Dipartimento di Matematica,
Universit\`a degli Studi ``Roma Tre", 00146 Roma, Italy. E-mail: esposito@mat.uniroma3.it. Research supported by M.U.R.S.T., project ``Variational methods and
nonlinear differential equations".}\qquad
Nassif Ghoussoub\footnote{Department of Mathematics, University of British Columbia, Vancouver, B.C. Canada V6T 1Z2. E-mail: nassif@math.ubc.ca. Research partially supported by the Natural Science and Engineering Research Council of Canada.}\qquad
Amir Moradifam\footnote{Department of Mathematics, University of British Columbia, Vancouver, B.C. Canada V6T 1Z2. E-mail: a.moradi@math.ubc.ca.}
}}
\date{\today}

\smallbreak
\maketitle

\begin{abstract} We study the regularity of the extremal solution of the semilinear biharmonic equation $\bi u=\f{\lambda}{(1-u)^2}$, which models a simple Micro-Electromechanical System (MEMS) device on  a ball $B\subset\IR^N$, under Dirichlet boundary conditions
$u=\partial_\nu u=0$ on $\partial B$.  We complete here the results of F.H. Lin and Y.S. Yang \cite{LY}  regarding the identification of a ``pull-in voltage"  $\la^*>0$ such that  a  stable classical solution $u_\la$ with $0<u_\la<1$ exists for $\la\in (0,\la^*)$, while there is none of any kind when $\la>\la^*$. Our main result asserts  that  the extremal solution $u_{\lambda^*}$ is regular $(\sup_B u_{\lambda^*}  <1)$ provided $ N \le 8$ while  $u_{\lambda^*} $ is singular ($\sup_B u_{\lambda^*} =1$) for $N \ge 9$, in which case $1-C_0|x|^{4/3}\leq u_{\lambda^*} (x) \leq  1-|x|^{4/3}$ on the unit ball, where $ C_0:= \left( \frac{\lambda^*}{\overline{\lambda}}\right)^\frac{1}{3}$ and $ \bar{\lambda}:= \frac{8}{9} (N-\frac{2}{3}) (N- \frac{8}{3})$.
\end{abstract}

\section{Introduction}
The following model has been proposed  for the  description of  the steady-state of a simple Electrostatic MEMS device:
\begin{equation}
\label{MEMS.1}
 \arraycolsep=1.5pt
\left\{ \begin{array}{ll}
\alpha \Delta^2 u = \left( \beta \int_\Omega | \nabla u|^2 dx + \gamma \right) \Delta u + \frac{ \lambda f(x)}{ (1-u)^2 \left( 1 + \chi \int_\Omega \frac{dx}{(1-u)^2} \right)}  &\quad \hbox{in }\Omega \\
0<u<1      &\quad \hbox{in } \Omega \\
u=\alpha \partial_\nu u =0  &\quad \hbox{on } \partial \Omega , \quad
\end{array} \right. \end{equation}
 where $ \alpha, \beta, \gamma, \chi \ge 0$, $ f \in C( \overline{\Omega},[0,1])$ are fixed, $ \Omega$ is a bounded domain in $ \IR^N$ and $ \lambda \ge 0$ is a varying parameter (see for example Bernstein and Pelesko \cite{BP}).  The function $ u(x)$ denotes the height above a point $ x \in \Omega\subset \IR^N$ of a dielectric membrane clamped on $ \pOm$,
once it deflects torwards a ground plate fixed at height $z=1$, whenever a positive voltage -- proportional to $\lambda$ -- is applied.

\medskip \noindent In studying this problem,
one typically makes various simplifying assumptions on the parameters $ \alpha, \beta, \gamma, \chi$, and the  first approximation of (\ref{MEMS.1}) that has been studied extensively so far  is the equation
\begin{eqnarray*}
\hskip 150pt
\left\{ \begin{array}{ll}
-\Delta u=   \lambda  \frac{f(x)}{(1-u)^2} &\text{in } \Omega\\
\hfill 0<u<1 \quad \quad &\text{in } \Omega\hskip 150pt (S)_{ \lambda,f} \\
\hfill u=0 \quad \quad \quad &\text{on }\partial \Omega,
\end{array} \right.
\end{eqnarray*}
  where we have set $ \alpha = \beta = \chi=0$ and $ \gamma=1$ (see for example \cite{Esp,EGG,GG} and the monograph \cite{EGG.book}) .   This simple model,  which lends itself to the vast literature on second order semilinear eigenvalue problems, is already a rich source of interesting mathematical problems. The case when the ``permittivity profile" $f$ is constant  ($f=1$) on a general domain  was studied in \cite{MP}, following the pioneering work of Joseph and Lundgren \cite{JL} who had considered the radially symmetric case. The case for a  non constant permittivity profile $ f$ was advocated
by Pelesko \cite{P}, taken up by \cite{GPW}, and studied in depth in \cite{Esp,EGG,GG}.  The starting point of the analysis is the existence of a pull-in voltage $\lambda^*(\Omega, f)$, defined as
$$ \lambda^*(\Omega, f):= \sup \Big\{ \lambda >0: \hbox{there exists a classical solution of } (S)_{\lambda, f} \Big\}.$$
It is then shown  that for every $ 0 < \lambda < \lambda^*$,  there exists a smooth minimal (smallest) solution of $(S)_{\lambda, f}$,
while for $ \lambda > \lambda^*$ there is no solution even in a weak sense.
Moreover, the branch $ \lambda \mapsto u_\lambda(x)$ is increasing for each $ x \in \Omega$,  and therefore the function
$u^*(x):= \lim_{\lambda \nearrow \lambda^*} u_\lambda(x)$ can be considered as a generalized solution that corresponds to the
pull-in voltage $\lambda^*$.  Now the issue of the regularity of this extremal solution -- which, by elliptic regularity theory, is equivalent to
whether $ \sup_\Omega u^*<1$ -- is  an important question for many reasons, not the least of which being the fact that it decides whether the set of solutions stops there, or whether a new branch of solutions emanates from a bifurcation state $(u^*,\lambda^*)$.  This issue turned out to depend closely on the dimension and on the permittivity profile $f$. Indeed, it was shown in \cite{GG} that $u^*$ is regular in dimensions $1\leq N\leq 7$, while it is not necessarily the case for $N\geq 8$. In other words, the dimension $N=7$ is critical for equation $(S)_\lambda$ (when $f=1$, we simplify the notation $(S)_{\lambda,1}$ into $(S)_\lambda$).  On the other hand, it is shown in \cite{EGG} that the regularity of $u^*$ can be restored in any dimension,  provided we allow for a power law profile $|x|^\eta$ with $\eta$ large enough.

\medskip \noindent The case where $ \beta = \gamma = \chi=0$ (and $ \alpha=1$) in the above model, that is when we are dealing with the following fourth order analog of $(S)_\lambda$
\begin{eqnarray*}
\hskip 150pt
\left\{ \begin{array}{ll}
\bi u=   \frac{\lambda}{(1-u)^2} &\text{in } \Omega\\
0<u<1 &\text{in } \Omega \hskip 150pt (P)_\lambda \\
u=\partial_\nu u=0 &\text{on }\partial \Omega,
\end{array} \right.
\end{eqnarray*}
was also considered by \cite{CDG,LY} but with limited success. One of the reasons is the lack of  a ``maximum principle" which plays such a crucial role in developing the theory for the Laplacian. Indeed, it is a well known fact that
    such a principle does not normally hold for general domains $ \Omega$ (at least for the clamped boundary conditions $ u = \partial_\nu u =0$ on $ \pOm$) unless one restricts attention to the unit ball $ \Omega =B$   in $ \IR^N$,  where one can exploit  a positivity preserving property of $\bi$ due to T. Boggio \cite{Boggio}.
  This is precisely what was done in the references mentioned above, where a  theory of the minimal branch associated with $(P)_\lambda$ is developed along the same lines as for $(S)_\lambda$. The second obstacle is the well-known difficulty of extracting energy estimates for solutions of  fourth order problems from their stability properties.  This means that the methods used to analyze the regularity of the extremal solution for $(S)_\lambda$ could not carry to the corresponding problem for $(P)_\lambda$.

\medskip \noindent This is the question we address in this paper as we eventually show the following result.

\begin{thm}  The unique extremal solution $u^*$ for $(P)_{\lambda^*}$ in $B$ is regular in dimension $1\leq  N \le 8$, while  it is singular (i.e, $ \sup_B u^*=1$)  for $ N \ge 9$.
\end{thm}
\noindent In other words, the critical dimension for $(P)_\lambda$ in $B$ is $N=8$, as opposed to being equal to $7$ in $(S)_\lambda$.  We add that our methods are heavily inspired by the recent paper of Davila et al. \cite{DDGM} where it is shown that $N=12$ is the critical dimension for the fourth order nonlinear eigenvalue problem
$$\left\{\begin{array}{ll}
\bi u=   \lambda e^u &\text{in } B\\
u=\partial_\nu u=0 &\text{on }\partial B,
\end{array} \right. $$
while the critical dimension for its second order counterpart (i.e.,  the Gelfand problem) is $N=9$. There is,  however, one major difference between our approach and the one used by Davila et al. \cite{DDGM}. It is related to the most delicate dimensions -- just above the critical one -- where they use a computer assisted proof  to establish the singularity of the extremal solution, while our method is more analytical and relies on improved and non standard Hardy-Rellich inequalities  recently established by Ghoussoub-Moradifam \cite{GM} (See Appendix).

\medskip \noindent Throughout this  paper, we will always consider problem $(P)_\lambda$ on the unit ball $B$. We start by recalling some of the results from \cite{CDG} concerning $(P)_\lambda$, that will be needed in the sequel.  We  define
$$\lambda^*:= \sup\Big\{ \lambda> 0: \hbox{ there exists a classical solution of }(P)_\lambda \Big\},$$
and note that we are not restricting our attention to radial solutions. We will deal also with weak solutions:
\begin{dfn} We say that $u$ is a weak solution of $(P)_\lambda $ if $ 0 \le u \le 1$ a.e. in $B$, $ \frac{1}{(1-u)^2} \in L^1(B)$ and
\[ \int_B u \Delta^2 \phi = \lambda \int_B \frac{\phi}{(1-u)^2}, \qquad \forall \phi \in C^4(\bar B) \cap H_0^2(B).\]
We say that $ u $ is a weak super-solution (resp. weak sub-solution) of $(P)_\lambda$ if the equality is replaced with the inequality $ \ge $ (resp. $ \le $) for all $\phi \in C^4(\bar B) \cap H_0^2(B)$ with $\phi \ge 0$.
\end{dfn}
\noindent We also introduce notions of regularity and stability.
\begin{dfn} Say that a weak solution $u$ of $(P)_\lambda $ is regular (resp. singular) if $\|u\|_\infty<1$ (resp. $=$) and stable (resp. semi-stable) if
$$\mu_1(u)=\inf \left\{ \int_B ( \Delta \phi)^2 -2 \lambda \int_B \frac{ \phi^2}{(1-u)^3}: \phi \in H_0^2(B), \| \phi \|_{L^2}=1 \right\}$$
is positive (resp. non-negative).
\end{dfn}
\noindent The following extension of Boggio's principle will be frequently used in the sequel (see \cite[Lemma 16]{AGGM} and \cite[Lemma 2.4]{DDGM}):
\begin{lemma}[Boggio's Principle]\label{boggio}
Let $u\in L^1(B)$.  Then $u\geq 0$ a.e. in $B$, provided one of the following conditions hold:
\begin{enumerate}

\item $u\in C^4(\overline{B})$,  $\bi u\geq 0$ on $B$, and  $u=\frac{\partial u}{\partial n}= 0$ on $\partial B$.
\item $\int_{B} u\bi\phi\,dx\geq 0$ for all $0\leq \phi \in C^4(\overline{B})\cap H_0^2(B)$.
\item $u\in H^2(B)$, $u=0$ and $\frac{\partial u}{\partial n} \leq 0$ on $\partial B$, and
 $\int_{B} \Delta u \Delta \phi \geq 0$  for all $0\leq \phi \in H^2_0(B)$.
\end{enumerate}
Moreover, either $u\equiv 0$ or $u>0$ a.e. in $B$.
\end{lemma}

\noindent The following  theorem summarizes the main results in \cite{CDG} that will be needed in the sequel:
\begin{thm}\label{CdG} The following assertions hold:
\begin{enumerate}
\item For each $ 0 < \lambda < \lambda^*$ there exists a classical minimal solution $u_\lambda$ of $ (P)_\lambda$. Moreover $ u_\lambda $ is radial and radially decreasing.
\item  For $ \lambda > \lambda^*$,  there are no weak solutions of $(P)_\lambda$.
\item For each $ x \in B$ the map $ \lambda \mapsto u_\lambda(x)$ is strictly increasing on $ (0,\lambda^*)$.
\item The pull-in voltage $ \lambda^*$ satisfies the following bounds:
\[ \max \left\{ \frac{ 32(10N-N^2-12)}{27}, \frac{8}{9} (N-\frac{2}{3}) (N- \frac{8}{3})\right\} \le \lambda^* \le \frac{4 \nu_1}{27}\]
where $ \nu_1$ denotes the first eigenvalue of $ \Delta^2 $ in $H_0^2(B)$.
\item For each $ 0 < \lambda < \lambda^*$, $u_\lambda$ is a stable solution (i.e., $ \mu_1(u_\lambda)>0$).
\end{enumerate}
\end{thm}
\noindent Using the stability of $ u_\lambda $,  it can be shown that  $ u_\lambda $ is uniformly bounded in $H_0^2(B)$ and that $ \frac{1}{1-u_\lambda} $ is uniformly bounded in $L^3(B)$.  Since now $ \lambda \mapsto u_\lambda(x)$ is increasing, the function $ u^*(x):= \lim_{ \lambda \nearrow \lambda^*} u_\lambda(x)$ is well defined (in the pointwise sense),  $ u^* \in H_0^2(B)$, $ \frac{1}{1-u^*} \in L^3(B)$ and $ u^*$ is a weak solution of $ (P)_{\lambda^*}$.   Moreover $ u^*$ is the unique weak solution of $(P)_{\lambda^*}$.

\medskip \noindent The second result we list from \cite{CDG} is critical in identifying the extremal solution.
\begin{thm}  If $ u \in H_0^2(B)$ is a singular weak solution of $(P)_\lambda$, then $ u $ is  semi-stable if and only if $ (u, \lambda) =(u^*,\lambda^*)$.
\end{thm}

\section{The effect of boundary conditions on the pull-in voltage}

As in \cite{DDGM}, we are led to examine problem $(P)_\lambda$ with non-homogeneous boundary conditions such as
$$\hskip 150pt
\left\{\begin{array}{ll}
\Delta^2 u= \frac{ \lambda}{(1-u)^2} &\hbox{in } B \\
\alpha<u<1 &\hbox{in }B \hskip 150 pt (P)_{\lambda, \alpha, \beta}\\
u= \alpha\:,\:\:\partial_\nu u = \beta &\hbox{on } \partial B, \end{array}\right.
$$
where $\alpha, \beta$ are given.

\medskip \noindent Notice first that some restrictions on $ \alpha $ and $ \beta$ are necessary. Indeed, letting $\Phi(x):=( \alpha - \frac{\beta}{2} ) + \frac{\beta}{2} |x|^2 $ denote the unique solution of
\begin{equation} \label{Phi} \left\{ \begin{array}{ll}
\Delta^2 \Phi = 0 &\hbox{in } B \\
\Phi = \alpha\:,\:\: \partial_\nu \Phi = \beta&\hbox{on }\partial B,
\end{array}\right.
\end{equation}
we infer immediately from Lemma \ref{boggio} that the function $u-\Phi$ is positive in $B$,  which yields to
$$\sup_B \Phi<\sup_B u\leq 1.$$
To insure that $\Phi$ is a classical sub-solution of $(P)_{\lambda,\alpha,\beta}$, we impose $\alpha \not= 1$ and $\beta \leq 0$, and condition $\displaystyle \sup_B \Phi<1$ rewrites as
$\alpha-\frac{\beta}{2} < 1$. We will then say that the pair $ (\alpha, \beta)$ is {\it admissible} if $\beta \leq 0$, and $\alpha-\frac{\beta}{2} < 1$.

\medskip \noindent This section will be devoted to obtaining  results for $ (P)_{ \lambda, \alpha ,\beta}$ when $ (\alpha, \beta)$ is an admissible pair, which are  analogous to those  for $ (P)_\lambda$.  To cut down on notation, we shall  sometimes drop $ \alpha $ and $ \beta$ from our expressions whenever such an emphasis is not needed.  For example in this section $ u_\lambda $ and $ u^*$ will denote the minimal and extremal solution of $ (P)_{\lambda, \alpha , \beta}$.

\medskip \noindent We now introduce a notion of weak solution for $(P)_{\lambda,\alpha,\beta}$.
\begin{dfn} We say that $u$ is a weak solution of $(P)_{\lambda,\alpha,\beta}$ if $\alpha \leq u \le 1$ a.e. in $B$, $ \frac{1}{(1-u)^2} \in L^1(B)$ and if
\[ \int_B (u-\Phi) \Delta^2 \phi = \lambda \int_B \frac{\phi}{(1-u)^2}, \qquad \forall \phi \in C^4(\bar B) \cap H_0^2(B),\]
where $\Phi$ is given in (\ref{Phi}). We say that $ u $ is a weak super-solution (resp. weak sub-solution) of $(P)_{\lambda,\alpha,\beta}$ if the equality is replaced with the inequality $ \ge $ (resp. $ \le $) for $ \phi \ge 0$.
\end{dfn}
\noindent We now define as before
\[ \lambda^*:= \sup \{ \lambda > 0: (P)_{\lambda, \alpha, \beta} \; \mbox{ has a classical solution} \}\]
and
\[\lambda_*:= \sup \{ \lambda > 0: (P)_{\lambda, \alpha, \beta} \; \mbox{ has a weak solution} \}.\]
Observe that by the Implicit Function Theorem,  one can always solve $(P)_{\lambda,\alpha,\beta}$ for small $\lambda$'s. Therefore, $\lambda^*$ (and also $\lambda_*$) is well defined.

\medskip \noindent Let now $U$ be a weak super-solution of $(P)_{\lambda,\alpha,\beta}$. Recall the following standard existence result.
\begin{thm} [\cite{AGGM}] \label{exist} For every $0\leq f \in L^1(B)$,  there exists a unique $0\leq u \in L^1(B)$ which satisfies
$$\int_B u \Delta^2 \phi=\int_B f \phi$$
for all $\phi \in C^4(\bar B) \cap H_0^2(B)$.
\end{thm}
\noindent We can now introduce the following ``weak iterative scheme": Start with $u_0=U$ and (inductively) let $u_n$, $n \geq 1$, be the solution of
$$\int_B (u_n-\Phi) \Delta^2 \phi=\lambda \int_B  \frac{\phi}{(1-u_{n-1})^2}\qquad\:\forall \: \phi \in C^4(\bar B) \cap H_0^2(B)$$
given by Theorem \ref{exist}. Since $0$ is a sub-solution of $(P)_{\lambda,\alpha,\beta}$, one can easily show inductively by using Lemma \ref{boggio} that $\alpha \leq u_{n+1}\leq u_n \leq U$ for every $n \geq 0$. Since
$$(1-u_n)^{-2}\leq (1-U)^{-2} \in L^1(B),$$
we get by Lebesgue Theorem, that the function $u=\displaystyle \lim_{n \to +\infty} u_n$ is a weak solution of $(P)_{\lambda,\alpha,\beta}$ such that $\alpha \leq u\leq U$. In other words,  the following result holds.

\begin{prop} \label{super} Assume the existence of a weak super-solution $U$ of $(P)_{\lambda,\alpha,\beta}$. Then there exists a weak solution $u$ of $(P)_{\lambda,\alpha,\beta}$ so that $\alpha \leq u \leq U$ a.e. in $B$.
\end{prop}
\noindent In particular,  we can find a weak solution of $(P)_{\lambda,\alpha,\beta}$ for every $\lambda \in (0,\lambda_*)$. Now we show that this is still true for regular weak solutions.

\begin{prop} \label{cch} Let $ (\alpha, \beta)$ be an admissible pair and let $u$ be a weak solution of $(P)_{\lambda,\alpha,\beta}$. Then for every $0<\mu<\lambda$, there is a regular solution for $(P)_{\mu,\alpha,\beta}$.
\end{prop}
\begin{proof} Let $ \E\in (0,1)$ be given and let  $ \bar u=(1-\E)u+\E \Phi$, where $\Phi$ is given in (\ref{Phi}). We have that
$$\sup_B \bar u\leq (1-\E)+\E \sup_B \Phi<1\:,\quad \inf_B \bar u\geq (1-\E)\alpha +\E \inf_B \Phi=\alpha,$$
and for every $0\leq \phi \in C^4(\bar B) \cap H_0^2(B)$ there holds:
\begin{eqnarray*}
\int_B (\bar u-\Phi) \Delta^2 \phi &=& (1-\E) \int_B (u-\Phi)\Delta^2 \phi
= (1-\E)\lambda \int_B \frac{\phi}{(1-u)^2}\\
&=& (1-\E)^3 \lambda \int_B \frac{\phi}{(1-\bar u+\E (\Phi-1))^2} \geq
(1-\E)^3 \lambda \int_B \frac{\phi}{(1-\bar u)^2}.
\end{eqnarray*}
Note that $ 0 \le (1-\E)(1-u)=1 - \bar{u}+\E (\Phi -1) <1-\bar u$. So $ \bar{u}$ is a weak super-solution of $ (P)_{ (1-\E)^3 \lambda, \alpha , \beta}$ satisfying $\displaystyle \sup_B \bar u<1$. From Proposition \ref{super} we get the existence of a weak solution $w$ of $ (P)_{ (1-\E)^3 \lambda, \alpha , \beta}$ so that $\alpha \leq w\leq \bar u$. In particular, $\displaystyle \sup_B w<1$ and $w$ is a regular weak solution. Since $\E \in (0,1)$ is arbitrarily chosen, the proof is complete.
\end{proof}
\noindent Proposition \ref{cch} implies in particular the existence of a regular weak solution $U_\lambda$ for every $\lambda \in (0,\lambda_*)$. Introduce now a ``classical" iterative scheme: $u_0=0$ and (inductively) $u_n=v_n +\Phi$, $n \geq 1$, where $v_n \in H_0^2(B)$ is the (radial) solution of
\begin{equation} \label{pranzo}
\Delta^2 v_n=\Delta^2(u_n-\Phi)= \frac{\lambda}{(1-u_{n-1})^2} \qquad\hbox{in }B.
\end{equation}
Since $v_n \in H_0^2(B)$, $u_n$ is also a weak solution of (\ref{pranzo}), and by Lemma \ref{boggio} we know that $\alpha \leq u_n\leq u_{n+1} \leq U_\lambda$ for every $n \geq 0$. Since $\displaystyle \sup_B u_n \leq \displaystyle \sup_B U_\lambda<1$ for $n\geq 0$, we get that $(1-u_{n-1})^{-2} \in L^2(B)$ and the existence of $v_n$ is guaranteed. Since $v_n$ is easily seen to be uniformly bounded in $H_0^2(B)$, we have that $u_\lambda:=\displaystyle \lim_{n \to +\infty}u_n$ does hold pointwise and weakly in $H^2(B)$. By Lebesgue Theorem, we have that $u_\lambda$ is a radial weak solution of $(P)_{\lambda,\alpha,\beta}$ so that $\displaystyle \sup_B u_\lambda\leq \displaystyle \sup_B U_\lambda<1$. By elliptic regularity theory \cite{ADN} $u_\lambda \in C^\infty(\bar B)$ and $u_\lambda-\Phi=\partial_\nu(u_\lambda-\Phi)=0$ on $\partial B$. So we can integrate by parts to get
$$\int_B \Delta^2 u_\lambda \phi =\int_B \Delta^2(u_\lambda-\Phi) \phi=\int_B (u_\lambda-\Phi)\Delta^2 \phi=\lambda \int_B \frac{\phi}{(1-u_\lambda)^2}$$
for every $\phi \in C^4(\bar B) \cap H_0^2(B)$. Hence, $u_\lambda$ is a radial classical solution of $(P)_{\lambda,\alpha,\beta}$ showing that $\lambda^*=\lambda_*$. Moreover, since $\Phi$ and $v_\lambda:=u_\lambda-\Phi$ are radially decreasing in view of \cite{Sor}, we get that $u_\lambda$ is radially decreasing too. Since the argument above shows that $u_\lambda<U$ for any other classical solution $U$ of $(P)_{\mu,\alpha,\beta}$ with $\mu \geq \lambda$, we have that $u_\lambda$ is exactly the minimal solution and $u_\lambda$ is strictly increasing as $\lambda \uparrow \lambda^*$. In particular, we can define $ u^*$ in the usual way: $ u^*(x)= \displaystyle \lim_{\lambda \nearrow \lambda^*} u_\lambda(x)$.

\medskip \noindent Finally, we show the finiteness of the pull-in voltage.
\begin{lemma} If $ (\alpha, \beta)$ is an admissible pair, then $\lambda^*(\alpha, \beta) <+\infty$.
\end{lemma}
\begin{proof} Let $u$ be a classical solution of $ (P)_{\lambda, \alpha, \beta}$  and let $ (\psi, \nu_1)$ denote the first eigenpair of $ \Delta^2$ in $H_0^2(B)$ with $ \psi >0$.  Now, let $ C $ be such that
\[ \int_{\partial B} (\beta \Delta \psi - \alpha \partial_\nu \Delta \psi) = C \int_B \psi. \]
Multiplying $ (P)_{\lambda,\alpha,\beta}$ by $ \psi$ and then integrating by parts one arrives at
\[ \int_B \left( \frac{ \lambda}{(1-u)^2} - \nu_1 u -C \right) \psi =0. \]
Since $ \psi>0$ there must exist a point $\bar x \in B$ where
$\frac{ \lambda}{(1-u(\bar x))^2} - \nu_1 u(\bar x) -C \le 0.$
Since $\alpha<u(\bar x)<1$, one can conclude that
$ \lambda \le \sup_{\alpha< u <1} ( \nu_1 u +C)(1-u)^2$,
which shows that $ \lambda^*<+\infty$.
\end{proof}
\noindent The following summarizes what we have shown so far.

\begin{thm} If $(\alpha,\beta)$ is an admissible pair, then $\lambda^*:=\lambda^*(\alpha, \beta) \in (0,+\infty)$ and the following hold:
\begin{enumerate}
\item For each $ 0 < \lambda < \lambda^*$ there exists a classical, minimal solution $u_\lambda$ of $(P)_{\lambda,\alpha,\beta}$.  Moreover $ u_\lambda $ is radial and radially decreasing.
\item For each $ x \in B$ the map $ \lambda \mapsto u_\lambda(x)$ is strictly increasing on $ (0,\lambda^*)$.
\item For $ \lambda > \lambda^*$ there are no weak solutions of $(P)_{\lambda,\alpha,\beta}$.
\end{enumerate}
\label{quasi} \end{thm}

\subsection{Stability of the minimal branch of solutions}
This section is devoted to the proof of the following stability result for minimal solutions.  We shall need yet another notion of $H^2(B)-$weak solutions, which is an intermediate class between classical and weak solutions.

\begin{dfn} We say that $u$ is a $H^2(B)-$weak solution of $(P)_{\lambda,\alpha,\beta}$ if $u -\Phi \in H_0^2(B)$, $ \alpha \le u \le 1$ a.e. in $B$, $ \frac{1}{(1-u)^2} \in L^1(B)$ and if
\[ \int_B \Delta u \Delta \phi = \lambda \int_B \frac{\phi}{(1-u)^2}, \qquad \forall \phi \in C^4(\bar B) \cap H_0^2(B),\]
where $\Phi$ is given in (\ref{Phi}). We say that $u$ is a $H^2(B)-$weak super-solution (resp. $H^2(B)-$weak sub-solution) of $(P)_{\lambda,\alpha,\beta}$ if for $ \phi \ge 0$ the equality is replaced with $ \ge$ (resp. $ \le $) and $u\geq \alpha$ (resp. $\leq$), $\partial_\nu u \leq \beta$ (resp. $\geq$) on $\partial B$.
\end{dfn}

\begin{thm} \label{stable} Suppose $ (\alpha,\beta)$ is an admissible pair.
\begin{enumerate}
\item  The minimal solution $ u_\lambda $ is then stable and is the unique semi-stable $H^2(B)-$weak solution of $(P)_{\lambda,\alpha,\beta}$.
\item The function $ u^*:= \displaystyle \lim_{\lambda \nearrow \lambda^*} u_\lambda$ is a well-defined semi-stable $H^2(B)-$weak solution of  $(P)_{\lambda^*,\alpha,  \beta}$.
\item When $u^*$ is classical solution, then $\mu_1(u^*)=0$ and $u^*$ is the unique $H^2(B)-$weak solution of $(P)_{\lambda^*,\alpha,\beta}$.
\item If $ v$ is a singular, semi-stable $H^2(B)-$weak solution of $ (P)_{ \lambda, \alpha, \beta}$, then $ v=u^*$ and $ \lambda = \lambda^*$
\end{enumerate}
\end{thm}

\noindent The crucial tool is a comparison result which is valid exactly in this class of solutions.
\begin{lemma} \label{shi} Let $ (\alpha, \beta)$ be an admissible pair and $u$ be a semi-stable $H^2(B)-$weak solution of $(P)_{\lambda, \alpha, \beta}$. Assume $ U $ is a $H^2(B)-$weak super-solution of $(P)_{\lambda, \alpha, \beta}$ so that $U-\Phi \in H_0^2(B)$. Then
\begin{enumerate}
\item $ u \le U$ a.e. in $B$;
\item If $ u$ is a classical solution and $ \mu_1(u)=0$ then $ U=u$.
\end{enumerate}
\end{lemma}
\begin{proof}  (i)  Define $ w:= u-U$.   Then by the Moreau decomposition \cite{M} for the biharmonic operator,  there exist $ w_1,w_2 \in H_0^2(B)$, with $ w=w_1 + w_2$, $ w_1 \ge 0$ a.e., $\Delta^2 w_2 \le 0 $ in the $H^2(B)-$weak sense and $\int_B \Delta w_1 \Delta w_2=0$.  By Lemma \ref{boggio}, we have that $w_2 \le 0$ a.e. in $B$.\\
Given now $ 0 \le \phi \in C_c^\infty(B)$, we have that
\[ \int_B \Delta w \Delta \phi \leq  \lambda \int_B (f(u) - f(U)) \phi, \]
where $ f(u):= (1-u)^{-2}$.  Since $ u$ is semi-stable, one has
\begin{eqnarray*}
\lambda \int_B f'(u) w_1^2 \le  \int_B (\Delta w_1)^2
= \int_B \Delta w \Delta w_1 \le  \lambda \int_B ( f(u) - f(U)) w_1.
\end{eqnarray*}
Since $ w_1 \ge w$ one also has
\[ \int_B f'(u) w w_1 \le \int_B (f(u)-f(U)) w_1,\]
which once re-arranged gives
\[ \int_B \tilde{f} w_1 \geq 0,\]
where $ \tilde{f}(u)= f(u) - f(U) -f'(u)(u-U)$. The strict convexity of $f$ gives $ \tilde{f} \le 0$ and $ \tilde{f}< 0 $ whenever $u \not= U$. Since $w_1 \ge 0$ a.e. in $B$ one sees that $ w \le 0 $ a.e.  in $B$. The inequality $ u \le U$  a.e. in $B$ is then established.

\medskip \noindent (ii)  Since $u$ is a classical solution, it is easy to see that the infimum in $\mu_1(u)$ is attained at some $\phi$. The function $\phi$ is then the first eigenfunction of $\Delta^2-\frac{2\lambda}{(1-u)^3}$ in $H_0^2(B)$. Now we show that $ \phi$ is of fixed sign.  Using the above decomposition, one has $ \phi= \phi_1 + \phi_2$ where $ \phi_i \in H_0^2(B)$ for $i=1,2$, $ \phi_1 \ge 0$, $ \int_B \Delta \phi_1 \Delta \phi_2=0$ and $ \Delta^2 \phi_2 \le 0$ in the $H^2_0(B)-$weak sense. If $ \phi$ changes sign,  then $ \phi_1 \not\equiv 0$ and $ \phi_2 <0$ in $B$ (recall that either $\phi_2<0$ or $\phi_2=0$ a.e. in $B$). We can write now:
\begin{eqnarray*}
0 = \mu_1(u)
\le  \frac{ \int_B (\Delta (\phi_1 -\phi_2))^2 - \lambda f'(u) ( \phi_1 - \phi_2)^2}{ \int_B ( \phi_1 - \phi_2)^2} < \frac{ \int_B ( \Delta \phi)^2 - \lambda f'(u) \phi^2 }{ \int_B \phi^2} =\mu_1(u)
\end{eqnarray*}
in view of $\phi_1 \phi_2<-\phi_1\phi_2$ in a set of positive measure, leading to a contradiction.\\
So we can assume $ \phi \ge 0$, and by the Boggi's principle we have $\phi>0$ in $B$. For $ 0 \le t \le 1$ define
$$g(t)=\int_B \Delta \left[t U+(1-t)u \right] \Delta \phi
- \lambda \int_B f( tU+(1-t)u) \phi,$$
where $\phi$ is the above first eigenfunction.
Since $ f$ is convex one sees that
$$g(t)\geq \lambda \int_B \left[t f(U)+(1-t)f(u)-f(tU+(1-t)u)\right]\phi \geq 0$$  for every $t \geq 0$. Since $ g(0) =0$ and
$$ g'(0)= \int_B \Delta (U-u) \Delta \phi-\lambda f'(u)(U-u)\phi=0 ,$$
we get that
\[ g''(0)=- \lambda \int_B f''(u) (U-u)^2 \phi\geq 0.\]
Since $f''(u)\phi>0$ in $B$, we finally get that $ U=u$ a.e. in $B$.
\end{proof}

\noindent Based again on Lemma \ref{boggio}(3), we can show a more general version of the above Lemma \ref{shi}.
\begin{lemma} \label{poo} Let $ (\alpha,\beta)$ be an admissible pair and $\beta'\leq 0$. Let $u$ be a semi-stable $H^2(B)-$weak sub-solution of $(P)_{\lambda, \alpha,\beta}$ with $u=\alpha$, $\partial_\nu u=\beta' \geq \beta$ on $\partial B$. Assume that $U$ is a $H^2(B)-$weak super-solution of $(P)_{\lambda, \alpha ,\beta}$ with $U=\alpha$, $\partial_\nu U=\beta$ on $\partial B$. Then $ U \ge u$ a.e. in $B$.
\end{lemma}
\begin{proof}  Let $ \tilde{u} \in H_0^2(B)$ denote a weak solution to $ \Delta^2 \tilde{u}= \Delta^2 (u-U)$ in $B$. Since $\tilde{u}-u+U=0$ and $\partial_\nu(\tilde{u}-u+U)\leq 0$ on $\partial B$, by Lemma \ref{boggio} one has that $ \tilde{u} \ge u-U $ a.e. in $B$.  Again by the Moreau decomposition \cite{M},  we may write $\tilde u$ as $ \tilde{u} = w+v $, where $ w,v \in H_0^2(B)$, $ w \ge 0 $ a.e. in $B$, $ \Delta^2 v \le 0$ in a $H^2(B)-$weak sense and $\int_B \Delta w \Delta v=0$.   Then for $ 0 \le \phi \in C^4 (\bar B)\cap H_0^2(B)$ one has
\[ \int_B \Delta \tilde{u} \Delta \phi =\int_B \Delta(u-U) \Delta \phi \leq \lambda \int_B (f(u)- f( U)) \phi .\]
In particular, we have that
\[ \int_B \Delta \tilde{u} \Delta w \le \lambda \int_B ( f(u)-f(U)) w.\]
Since by semi-stability of $u$
\begin{eqnarray*}
\lambda \int_B f'(u) w^2\leq \int_B ( \Delta w)^2 = \int_B \Delta \tilde{u} \Delta w ,
\end{eqnarray*}
we get that
\[  \int_B f'(u) w^2 \le \int_B ( f(u)-f(U)) w.\]  By Lemma \ref{boggio} we have $v\leq 0$ and then $ w \ge \tilde{u} \ge u -U$ a.e. in $B$. So we see that
\[ 0 \le \int_B \left( f(u)-f(U)-f'(u)(u-U) \right) w.\]  The strict convexity of $ f$ implies as in Lemma \ref{shi} that $ U \ge u $ a.e. in $B$.
\end{proof}
\noindent We shall need the following a-priori estimates along the minimal branch $u_\lambda$.

\begin{lemma} \label{extremalsol} Let $ (\alpha, \beta)$ be an admissible pair. Then one has  \[ 2 \int_B \frac{( u_\lambda - \Phi)^2}{(1-u_\lambda)^3} \le \int_B \frac{ u_\lambda - \Phi}{(1-u_\lambda)^2},\]
where $ \Phi$ is given in (\ref{Phi}). In particular, there is a constant $C>0$ so that for every $\lambda \in (0,\lambda^*)$, we have
\begin{equation} \label{tardi}
\int_B (\Delta u_\lambda)^2+\int_B \frac{1}{(1-u_\lambda)^3} \leq C.
\end{equation}

\end{lemma}
\begin{proof} Testing $ (P)_{\lambda, \alpha , \beta}$ on $ u_\lambda - \Phi \in C^4(\bar B) \cap H^2_0(B)$, we see that
\begin{eqnarray*}
\lambda \int_B \frac{ u_\lambda - \Phi}{(1-u_\lambda)^2} = \int_B \Delta u_\lambda \Delta( u_\lambda - \Phi) =\int_B ( \Delta (u_\lambda - \Phi))^2
\ge 2 \lambda \int_B \frac{ (u_\lambda- \Phi)^2}{( 1-u_\lambda)^3}
\end{eqnarray*}
in view of $\Delta^2 \Phi=0$. In particular, for $\delta>0$ small we have that
\begin{eqnarray*}
\int_{\{|u_\lambda-\Phi| \geq \delta \}}\frac{1}{(1-u_\lambda)^3}&\leq &
\frac{1}{\delta^2} \int_{\{|u_\lambda-\Phi| \geq \delta \}}\frac{(u_\lambda-\Phi)^2}{(1-u_\lambda)^3}
\leq \frac{1}{\delta^2} \int_B \frac{1}{(1-u_\lambda)^2}\\
&\leq &\delta \int_{\{|u_\lambda-\Phi| \geq \delta \}}\frac{1}{(1-u_\lambda)^3}+C_\delta
\end{eqnarray*}
by means of Young's inequality.
Since for $\delta$ small,
$$\int_{\{|u_\lambda-\Phi| \leq \delta \}}\frac{1}{(1-u_\lambda)^3}\leq C'$$
for some $C'>0$, we  can deduce that for every $\lambda \in (0,\lambda^*)$,
$$\int_B \frac{1}{(1-u_\lambda)^3} \leq C$$
for some $C>0$.  By Young's and H\"older's inequalities, we now have
$$\int_B (\Delta u_\lambda)^2=\int_B \Delta u_\lambda \Delta \Phi+\lambda \int_B \frac{u_\lambda -\Phi}{(1-u_\lambda)^2}\leq \delta \int_B (\Delta u_\lambda)^2
+C_\delta+C \left(\int_B \frac{1}{(1-u_\lambda)^3} \right)^{\frac{2}{3}}$$
 and estimate (\ref{tardi}) is therefore established.\end{proof}

\medskip \noindent We are now ready to establish Theorem \ref{stable}.\\

\noindent {\bf Proof (of Theorem \ref{stable}):} (1)\, Since $\|u_\lambda\|_\infty <1$, the infimum defining $\mu_1(u_\lambda)$ is achieved at a first eigenfunction for every $\lambda \in (0,\lambda^*)$. Since $\lambda  \mapsto u_\lambda(x)$ is increasing for every $x \in B$, it is easily seen that $\lambda \mapsto \mu_1( u_\lambda)$ is an increasing, continuous function on $ (0, \lambda^*)$.  Define
\[ \lambda_{**}:= \sup\{ 0 <\lambda < \lambda^*: \: \mu_1( u_\lambda) >0 \} .\] We have that $ \lambda_{**}= \lambda^*$. Indeed, otherwise we would have that $ \mu_1(u_{ \lambda_{**}}) =0$, and for every $ \mu \in ( \lambda_{**}, \lambda^*)$ $ u_{\mu}$ would be a classical super-solution of $ (P)_{ \lambda_{**},\alpha, \beta}$. A contradiction arises since Lemma \ref{shi} implies $u_{\mu} = u_{\lambda_{**}}$.\\
Finally, Lemma \ref{shi} guarantees uniqueness in the class of semi-stable $H^2(B)-$weak solutions.\\
(2) \, By estimate (\ref{tardi}) it follows that $u_\lambda \to u^*$ in a pointwise sense and weakly in $H^2(B)$, and $ \frac{1}{1-u^*} \in L^3(B)$. In particular, $u^*$ is a $H^2(B)-$weak solution of $(P)_{ \lambda^*, \alpha ,\beta}$ which is also semi-stable as limiting function of the semi-stable solutions $\{u_\lambda\}$.\\
(3) Whenever $\|u^*\|_\infty<1$, the function $u^*$ is a classical solution, and by the Implicit Function Theorem  we have that $\mu_1(u^*)=0$ to prevent the continuation of the minimal branch beyond $\lambda^*$. By Lemma \ref{shi} $u^*$ is then the unique $H^2(B)-$weak solution of $(P)_{\lambda^*,\alpha,\beta}$. An alternative approach --which we do not pursue here-- based on the very definition of the extremal solution $u^*$ is available in \cite{CDG} when $\alpha=\beta=0$ (see also \cite{Mar}) to show that $u^*$ is the unique weak solution of $(P)_{\lambda^*}$, regardless of whether $u^*$ is regular or not.\\
(4) \, If $ \lambda < \lambda^*$, by uniqueness $v=u_\lambda $. So $v$ is not singular and a contradiction arises.

\medskip \noindent By Theorem \ref{quasi}(3) we have that $ \lambda = \lambda^*$. Since $ v $ is a semi-stable $H^2(B)-$weak solution of $ (P)_{ \lambda^*, \alpha, \beta}$ and $ u^*$ is a $H^2(B)-$weak super-solution of $ (P)_{\lambda^*, \alpha , \beta}$, we can apply Lemma \ref{shi} to get $ v \le u^*$ a.e. in $B$.  Since $u^*$ is a semi-stable solution too, we can reverse the roles of $ v$ and $ u^*$ in Lemma \ref{shi} to see that $ v \ge u^*$ a.e. in $B$. So equality $v=u^*$ holds and the proof is done.

\section{Regularity of the extremal solution  for $ 1 \le N \le 8$ }
We now return to the issue of the regularity of the extremal solution in problem $(P)_\lambda$. Unless stated otherwise, $ u_\lambda $ and $ u^*$ refer to the minimal and extremal solutions of $ (P)_\lambda$.   We shall show that the extremal solution $ u^*$ is regular provided $ 1 \le N \le 8$.
We first begin by showing that it is indeed the case  in small dimensions:
\begin{thm} $ u^*$ is regular in dimensions $ 1 \le N \le 4$.
\label{regular1} \end{thm}
\begin{proof} As already observed, estimate (\ref{tardi}) implies that $f(u^*)=(1-u^*)^{-2} \in L^{\frac{3}{2}}(B)$. Since $u^*$ is radial and radially decreasing, we need to show that $ u^*(0)<1$ to get the regularity of $ u^*$. The integrability of $f(u^*)$ along with elliptic regularity theory shows that $ u^* \in W^{4, \frac{3}{2}}(B)$. By the Sobolev imbedding Theorem we get that $u^*$ is a Lipschitz function in $B$.\\
Now suppose $ u^*(0)=1$ and $ 1 \le N \le 3$. Since
$$\frac{1}{1-u} \ge \frac{C}{|x|}\qquad \hbox{in }B$$
for some $ C>0$, one sees that
\[
+\infty = C^3 \int_B \frac{1}{|x|^3} \le \int_B \frac{1}{(1-u^*)^3} < +\infty.\]
A contradiction arises and hence $u^*$ is regular for $ 1 \le N \le 3$.\\
For $N=4$ we need to be more careful and observe that $u^* \in C^{1, \frac{1}{3}}(\bar B)$ by the Sobolev Imbedding Theorem. If $ u^*(0)=1$, then $ \nabla u^*(0)=0$ and
\[ \frac{1}{1-u^*} \ge \frac{C}{|x|^\frac{4}{3}} \qquad \hbox{in }B \]
for some $ C>0$. We now obtain a contradiction exactly as above.
\end{proof}
\noindent We now tackle the regularity of $ u^*$ for $ 5 \le N \le 8$. We start with the following crucial result:
\begin{thm} Let $ N \ge 5$ and $ (u^*, \lambda^*)$ be the extremal pair of $(P)_\lambda$. When $u^*$ is singular, then
\[ 1-u^*(x) \le C_0 |x|^\frac{4}{3} \qquad \hbox{in }B,\]
where $ C_0:= \left( \frac{\lambda^*}{\overline{\lambda}}\right)^\frac{1}{3}$ and $ \bar{\lambda}=\bar{\lambda}_N:= \frac{8}{9} (N-\frac{2}{3}) (N- \frac{8}{3})$.
\label{touchdown}\end{thm}
\begin{proof} First note that Theorem \ref{CdG}(4) gives the lower bound:
\begin{equation}\label{lowbound}
\lambda^* \geq \bar \lambda= \frac{8}{9} (N-\frac{2}{3}) (N- \frac{8}{3}).
\end{equation}
For $ \delta >0$, we define $ u_\delta(x):=1-C_\delta |x|^\frac{4}{3}$ with $ C_\delta:= \left( \frac{\lambda^*}{\bar \lambda}+\delta \right)^\frac{1}{3}>1$.   Since $N\geq 5$, we have that $ u_\delta \in H^2_{loc}(\IR^N)$, $ \frac{1}{1-u_\delta} \in L^3_{loc}(\IR^N)$ and $ u_\delta $ is a $H^2-$weak solution of
\[ \Delta^2 u_\delta = \frac{ \lambda^* + \delta \bar{ \lambda}}{ (1-u_\delta)^2} \qquad \mbox{ in } \IR^N.\]
We claim that $u_\delta \leq u^*$ in $B$, which will finish  the proof  by just letting $\delta \to 0$.

\medskip \noindent Assume by contradiction that the set $\Gamma:=\{ r \in (0,1):u_\delta(r) >u^*(r) \}$ is non-empty, and let $r_1=\displaystyle \sup \:\Gamma$.
Since
\[ u_\delta(1) = 1 - C_\delta<0=u^*(1),\]
we have that $0 < r_1 < 1$ and one infers that
\[ \alpha:= u^*(r_1)=u_\delta(r_1) \:, \quad \beta:=( u^*)'(r_1) \geq u_\delta'(r_1) .\]
Setting $u_{\delta,r_1}(r)=r_1^{-\frac{4}{3}}\left(u_\delta(r_1 r)-1 \right)  +1$, we easily see that $u_{\delta,r_1}$ is a $H^2(B)-$weak super-solution of $(P)_{\lambda^*+\delta {\bar \lambda}_N,\alpha',\beta'}$, where
$$\alpha':= r_1^{-\frac{4}{3}}( \alpha-1) +1\:,\quad \beta':=
r_1^{-\frac{1}{3}} \beta.$$

\medskip \noindent Similarly, let us define $u^*_{r_1}(r)= r_1^{-\frac{4}{3}}\left( u^*(r_1 r)-1\right) +1$. The dilation map
\begin{equation}
w \to w_{r_1}(r)=r_1^{-\frac{4}{3}}\left( w(r_1 r)-1\right) +1
\end{equation}
 is a correspondence between solutions of $(P)_{\lambda}$ on $B$ and of $(P)_{\lambda,1-r_1^{-\frac{4}{3}},0}$ on $B_{r_1^{-1}}$ which preserves the $H^2-$integrability. In particular, $(u^*_{r_1},\lambda^*)$ is  the extremal pair of $(P)_{\lambda,1-r_1^{-\frac{4}{3}},0}$ on $B_{r_1^{-1}}$ (defined in the obvious way).
Moreover, $u^*_{r_1}$ is a singular semi-stable $H^2(B)-$ weak solution of $(P)_{\lambda^* ,\alpha',\beta'}$.

\medskip \noindent Since $u^*$ is radially decreasing, we have that $ \beta' \le 0$. Define the function $w$ as $w(x):= ( \alpha' - \frac{\beta'}{2}) + \frac{ \beta'}{2} |x|^2 + \gamma(x) $, where $ \gamma $ is a solution of $ \Delta^2 \gamma= \lambda^* $ in $B$ with $ \gamma = \partial_\nu \gamma =0 $ on $ \partial B$.  Then $ w$ is a classical solution of
$$ \left\{ \begin{array}{ll}
\Delta^2 w = \lambda^* &\hbox{in } B \\
w = \alpha'\:, \quad \partial_\nu w = \beta' &\hbox{on } \partial B.
\end{array} \right.$$
Since $\frac{\lambda^*}{(1-u^*_{r_1})^2}\geq \lambda^*$, by Lemma \ref{boggio} we have $ u^*_{r_1} \ge w $ a.e. in $B$.  Since $ w(0) = \alpha' - \frac{\beta'}{2}+ \gamma(0)$ and $ \gamma(0)>0$, the bound $ u^*_{r_1} \le 1$ a.e. in $B$ yields to $ \alpha' - \frac{\beta'}{2}<1$. Namely, $ (\alpha', \beta')$ is an admissible pair and by Theorem
\ref{stable}(4) we get that $(u^*_{r_1},\lambda^*)$ coincides with the extremal pair of $(P)_{ \lambda, \alpha', \beta'}$ in $B$.

\medskip \noindent Since $(\alpha',\beta')$ is an admissible pair and $u_{\delta,r_1}$ is a $H^2(B)-$weak super-solution of $(P)_{\lambda^*+\delta {\bar \lambda}_N,\alpha',\beta'}$, by Proposition \ref{super} we get the existence of a weak solution of $(P)_{\lambda^*+\delta {\bar \lambda}_N,\alpha',\beta'}$. Since $\lambda^*+\delta {\bar \lambda}_N>\lambda^*$, we contradict the fact that $\lambda^*$ is the extremal parameter of $(P)_{\lambda,\alpha',\beta'}$.
\end{proof}
\noindent Thanks to this lower estimate on $u^*$, we get the following result.
\begin{thm} If $ 5 \le N \le 8$, then the extremal solution $u^*$ of $ (P)_\lambda$ is regular.
\label{regular2} \end{thm}
\begin{proof} Assume that $ u^*$ is singular. For  $ \E>0$ set  $\psi(x):= |x|^{ \frac{4-N}{2}+\E}$ and note that \[
 (\Delta \psi)^2 = (H_N +O( \E)) |x|^{-N+2\E}
 \]
where
\[
H_N:= \frac{N^2 (N-4)^2}{16}.
\]
Given $\eta \in C_0^\infty(B)$, and since $N\geq 5$, we can use the test function $\eta \psi \in H_0^2(B)$ into the stability inequality to obtain
\[
2 \lambda \int_B \frac{\psi^2}{(1-u^*)^3} \le \int_B (\Delta \psi)^2 +O(1),
\]
where $O(1)$ is a bounded function as $ \E \searrow 0$. By Theorem \ref{touchdown} we find that
\[   2 {\bar \lambda}_N \int_B \frac{\psi^2}{|x|^4}  \le \int_B (\Delta \psi)^2 +O(1),\]
and then
\[   2 {\bar \lambda}_N \int_B |x|^{-N+2\E} \le (H_N +O(\E)) \int_B |x|^{-N+2\E} +O(1).\] Computing the integrals one arrives at
\[ 2 {\bar \lambda}_N \le H_N +O(\E).\]
As $ \E \to 0$ finally we obtain $2 {\bar \lambda}_N \le H_N$. Graphing this relation one sees that $ N \ge 9$.
\end{proof}
\noindent We can now slightly improve the lower bound (\ref{lowbound}).
\begin{cor} \label{lambda.bar}  In any dimension $N\geq 1$, we have
\begin{equation}\label{lower}
\lambda^*>{\bar \lambda}_N=\frac{8}{9} (N-\frac{2}{3}) (N- \frac{8}{3}).
\end{equation}
\end{cor}
\begin{proof} The function $\bar{u}:=1-|x|^\frac{4}{3}$ is a $H^2(B)-$ weak solution of $(P)_{{\bar \lambda}_N,0,-\frac{4}{3}}$. If by contradiction $\lambda^*={\bar \lambda}_N$, then $\bar u$ is a $H^2(B)-$weak super-solution of $(P)_\lambda$ for every $\lambda \in (0,\lambda^*)$. By Lemma \ref{shi} we get that $u_\lambda \le \bar{u}$ for all $ \lambda < \lambda^*$, and then $u^*\le \bar{u}$ a.e. in $B$.

\medskip \noindent If $1\leq N \leq 8$, $u^*$ is then regular by Theorems \ref{regular1} and \ref{regular2}. By Theorem \ref{stable}(3) there holds $\mu_1(u^*)=0$. Lemma \ref{shi} then yields that $u^*=\bar u$, which is a contradiction since then $u^*$ will not satisfy the boundary conditions.

\medskip \noindent If now $N\geq 9$ and $ \bar{\lambda} = \lambda^*$, then $C_0=1$ in Theorem \ref{touchdown}, and we then have  $ u^* \geq \bar{u}$. It means again that $u^*=\bar u$, a contradiction that completes the proof.
\end{proof}

\section{The extremal solution is singular for $N \ge 9$}

We prove in this section that  the extremal solution is singular for $N \ge 9$. For that we have to distinguish between three different ranges for the dimension. For each range, we will need a suitable Hardy-Rellich type inequality that will be established in the appendix, by using the recent results of Ghoussoub-Moradifam \cite{GM}.  As in the previous section $(u^*,\lambda^*)$ denotes the extremal pair of $ (P)_\lambda$.
\begin{trivlist}

\item {\bf $\bullet$ Case  $N\geq 17$:}\, To establish the singularity of $u^*$ for these dimensions we shall need the following well known improved Hardy-Rellich inequality, which is valid for $N \ge 5$. There exists $C>0$, such that  for all $\phi \in H_0^2(B)$ 
\begin{equation}\label{HR1}
 \int_B (\Delta \phi)^2\, dx \ge \frac{ N^2(N-4)^2}{16} \int_B \frac{\phi^2}{|x|^4} \, dx + C \int_B \phi^2 \, dx.
 \end{equation}

\item {\bf $\bullet$ Case  $10\leq N\leq 16$:}\,  For this case,  we shall need the following inequality valid  for all $\phi \in H^2_{0}(B)$
\begin{eqnarray} \label{HR2}
\int_{B}(\Delta \phi)^2 &\geq&
\frac{(N-2)^2(N-4)^2}{16}\int_{B}\frac{\phi^2}{(|x|^2-|x|^{\frac{N}{2}+1})(|x|^2-|x|^{\frac{N}{2}})}\\
&&+\frac{(N-1)(N-4)^2}{4}
\int_{B}\frac{\phi^2}{|x|^2(|x|^2-|x|^{\frac{N}{2}})}.\nonumber
\end{eqnarray}

\item {\bf $\bullet$ Case  $N=9$:}\,  This case is the trickiest and  will require the following inequality  for all $\phi \in H^2_{0}(B)$
\begin{equation}\label{HR3}
\int_{B}(\Delta \phi)^2\geq \int_{B}Q(|x|)\left(
P(|x|)+\frac{N-1}{|x|^2}\right)\phi^2,
\end{equation}
where
\[
\hbox{$P(r)=\frac{\Delta_N \varphi}{\varphi}$ \quad  and \quad  $Q(r)=\frac{\Delta_{N-2}\psi}{\psi}$,}
\]
with $\varphi$ and $\psi$ being two appropriately chosen polynomials, namely
\[
\hbox{$\varphi(r):=r^{-\frac{N}{2}+1}+r-1.9$\quad  and \quad
$\psi(r):=r^{-\frac{N}{2}+2}+20r^{-1.69}+10r^{-1}+10r+7r^2-48$.}
\]
Recall that for a radial function $\varphi$, we set $\Delta_N \varphi(r)=\varphi''(r)+\frac{(N-1)}{r}\varphi'(r)$.
\end{trivlist}

We shall first show the following upper bound on $u^*$.
\begin{lemma} \label{blow} If $ N \ge 9$, then $ u^* \le 1 - |x|^\frac{4}{3}$ in $B$.
\end{lemma}
\begin{proof} Recall from Corollary \ref{lambda.bar} that $ \bar{\lambda}:=\frac{8}{9} (N-\frac{2}{3}) (N- \frac{8}{3})< \lambda^*$.
We now claim that $ u_\lambda \le \bar{u}$ for all $ \lambda \in ( \bar{\lambda}, \lambda^*)$.  Indeed,  fix such a $ \lambda $ and assume by contradiction that
\[ R_1:= \inf \{ 0 \le R \le 1: u_\lambda < \bar{u} \mbox{ in } (R,1) \}>0.\]
From the boundary conditions, one has that $ u_\lambda(r) < \bar{u}(r)$ as $ r\to 1^-$. Hence, $0<R_1<1$, $ \alpha:=u_\lambda(R_1)=\bar{u}(R_1)$ and $ \beta:=u_\lambda'(R_1) \le \bar{u}'(R_1)$.  Introduce,  as in the proof of Theorem \ref{touchdown}, the functions  $u_{\lambda, R_1}$ and $\bar u_{R_1}$.  We have that $u_{\lambda,R_1}$ is a classical super-solution of $(P)_{{\bar \lambda}_N,\alpha',\beta'}$, where
$$\alpha':= R_1^{-\frac{4}{3}}( \alpha-1) +1\:,\quad \beta':=
R_1^{-\frac{1}{3}} \beta.$$
Note that $\bar u_{R_1}$ is a $H^2(B)-$weak sub-solution of $(P)_{{\bar \lambda}_N,\alpha',\beta'}$ which is also semi-stable in view of the Hardy-Rellich inequality (\ref{HR1}) and the fact that
\[
2{\bar \lambda}_N \leq H_N:= \frac{ N^2(N-4)^2}{16}.
\]
By Lemma \ref{poo},  we deduce that $u_{\lambda,R_1}\geq \bar u_{R_1}$ in $B$. Note that, arguing as in the proof of Theorem \ref{touchdown}, $(\alpha',\beta')$ is an admissible pair.
We have therefore shown that $u_\lambda \geq \bar u$ in $B_{R_1}$ and a contradiction arises in view of the fact that $\displaystyle \lim_{x \to 0} \bar u(x)=1$ and $\|u_\lambda\|_\infty<1$. It follows that  $u_\lambda \leq \bar u$ in $B$ for every $\lambda \in ({\bar \lambda}_N, \lambda^*)$, and in particular $u^*\leq \bar u$ in $B$.
\end{proof}
The following lemma is the key for the proof of the singularity of $u^*$ in higher dimensions.

\begin{lemma}\label{sing-lem}  Let $N\geq 9$. Suppose there exist $\lambda'>0$, $\beta>0$ and a singular radial function $w \in
H^{2}(B)$with $\frac{1}{1-w}\in L^{\infty}_{loc}(\bar B\setminus
\{0\})$ such that
\begin{equation}\label{cond1}
\left\{\begin{array}{ll} \Delta^2 w\leq \frac{\lambda'}{(1-w)^2} &\hbox{for }0<r<1,\\
w(1)=0, \  w'(1)=0,& \end{array} \right.
\end{equation}
and
\begin{equation}\label{cond2}
2\beta \int_{B}\frac{ \phi^2}{(1-w)^3}\leq \int_{B}(\Delta
\phi)^2\ \ \hbox{for all} \ \ \phi \in H^2_{0}(B),
\end{equation}
\begin{enumerate}
\item If $\beta\geq \lambda'$, then
$\lambda^{*} \leq \lambda'$.
\item If either $\beta > \lambda'$ or if $\beta = \lambda'=\frac{H_N}{2}$, then the extremal solution $u^{*}$ is necessarily singular.
\end{enumerate}
\end{lemma}
\noindent {\bf Proof:} 1) First, note that (\ref{cond2}) and
$\frac{1}{1-w}\in L^{\infty}_{loc}(\bar B\setminus
\{0\})$ yield to $\frac{1}{(1-w)^2}\in L^1(B)$. By a density argument, (\ref{cond1}) implies now that $w$ is a $H^2(B)-$weak sub-solution of $(P)_{\lambda'}$ whenever $N\geq 4$. If now $\lambda' <\lambda^*$, then by Lemma \ref{poo} $w$ would necessarily be below the minimal solution $u_{\lambda'}$, which is a contradiction since $w$ is singular while $u_{\lambda'}$ is regular. \\
2)  Suppose first that {\bf $\beta = \lambda'=\frac{H_N}{2}$} and that $N\geq 9$. Since by part 1) we have $\lambda^*\leq \frac{H_N}{2}$, we get from Lemma \ref{blow} and the improved Hardy-Rellich inequality (\ref{HR1}) that  there exists $C>0$ so that for all $ \phi \in H^2_0(B)$
\begin{eqnarray*}
\int_B (\Delta \phi)^2 - 2 \lambda^* \int_B \frac{ \phi^2}{(1-u^*)^3} \ge \int_B (\Delta \phi)^2 - H_N \int_B \frac{\phi^2}{|x|^4} \ge C\int_B \phi^2.
\end{eqnarray*}
It follows that $ \mu_1(u^*)>0$ and $ u^*$ must therefore be singular since otherwise, one could use the Implicit Function Theorem  to continue the minimal branch beyond $ \lambda^*$. \\
Suppose now that  $\beta > \lambda'$, and let
$\frac{\lambda'}{\beta}<\gamma<1$ in such a way that
\begin{equation}
\alpha:=(\frac{\gamma\lambda^{*}}{\lambda'})^{1/3}<1.
\end{equation}
Setting $\bar{w}:=1-\alpha(1-w)$, we claim that
\begin{equation}\label{claim.000}
u^{*}\leq \bar{w} \ \ \hbox{in} \ \ B.
\end{equation}
Note that by the choice of $\alpha$ we have
$\alpha^3\lambda'<\lambda^{*}$, and therefore
to prove (\ref{claim.000}) it suffices to show that for
$\alpha^3\lambda'\leq \lambda<\lambda^{*}$,
we have
$u_{\lambda}\leq \bar{w} \ \ \hbox{in} \ \ B.$
Indeed, fix such $\lambda$ and note that
\begin{eqnarray*}
\Delta^2\bar{w}=\alpha \Delta^2 w \leq \frac{
\alpha\lambda'}{(1-w)^2}=\frac{\alpha^3\lambda'}{(1-\bar{w})^2}\leq \frac{\lambda}{(1-\bar{w})^2}.
\end{eqnarray*}
Assume that  $u_\lambda \leq \bar w$ does not hold in $B$, and
consider
\begin{equation*}
R_{1}:=\sup_{}\{0\leq R\leq 1\mid u_{\lambda}(R)>\bar{w}(R)\} >0.
\end{equation*}
 Since $\bar w(1)=1-\alpha>0=u_\lambda(1)$, we then have $R_{1}<1$,
$u_{\lambda}(R_{1})=\bar{w}(R_{1})$ and $(u_{\lambda})'(R_1) \leq
(\bar{w})'(R_1)$. Introduce,  as in the proof of Theorem
\ref{touchdown}, the functions  $u_{\lambda,R_1}$ and $\bar
w_{R_1}$.  We have that $u_{\lambda,R_1}$ is a classical solution of
$(P)_{\lambda,\alpha',\beta'}$, where
$$\alpha':= R_1^{-\frac{4}{3}}( u_\lambda(R_1)-1) +1\:,\quad \beta':=
R_1^{-\frac{1}{3}} (u_\lambda)'(R_1).$$
Since $\lambda<\lambda^{*}$ and then
\[\frac{2\lambda}{(1-\bar{w})^3}\leq \frac{2\lambda^{*}}{\alpha^3(1-w)^3}=\frac{2\lambda'}{ \gamma(1-w)^3}<\frac{2\beta}{(1-w)^3},\]
by (\ref{cond2}) $\bar w_{R_1}$ is a stable $H^2(B)-$weak
sub-solution of $(P)_{\lambda,\alpha',\beta'}$. By Lemma \ref{poo},
we deduce that $u_\lambda \geq \bar w$ in $B_{R_1}$ which is
impossible, since $\bar{w}$ is singular while $u_{\lambda}$ is
regular. Note that, arguing as in the proof of Theorem
\ref{touchdown}, $(\alpha',\beta')$ is an admissible pair. This
establishes claim (\ref{claim.000}) which, combined with the above
inequality, yields
\[
\frac{2\lambda^{*}}{(1-u^{*})^3}\leq \frac{2\lambda^{*}}{ \alpha^3(1-w)^3}<\frac{2\beta}{(1-w)^3},
\]
and therefore
\[\inf_{\phi \in H^2_{0}(B)}\frac{\int_{B}(\Delta \phi)^2-\frac{2\lambda^{*}\phi^2}{(1-u^{*})^3}}{\int_{B}\phi^2}>0.\]
It follows that again $ \mu_1(u^*)>0$ and $ u^*$ must be singular, since otherwise, one could use the Implicit Function Theorem  to continue the minimal branch beyond $ \lambda^*$.

Consider for any $m>0$ the following function:
\begin{equation}
 w_m:=1-\frac{3m}{3m-4}r^{4/3}+\frac{4}{3m-4} r^m,
\end{equation}
which satisfies the right boundary conditions: $w_m(1)=w_m'(1)=0$. We can now prove that the extremal solution is singular for $N\geq 9$.

\begin{thm} Let $N\geq 9$. The following upper bounds on $\lambda^*$ hold:
\begin{enumerate}
\item If $N \ge 31$, then Lemma \ref{sing-lem} holds with $w:=w_2$, $\lambda'=27 \bar \lambda_N$ and $\beta=\frac{H_N}{2}$, and therefore  $\lambda^*(N) \le 27 \bar \lambda_N$.
\item If $17 \le N \le 30$, then Lemma \ref{sing-lem} holds with $w:=w_3$, $\lambda'=\beta=\frac{H_N}{2}$, and therefore $\lambda^*(N) \le \frac{H_N}{2}$.
\item If $10\le N \le 16$, then Lemma \ref{sing-lem} holds with $w:=w_3$, $\lambda'_N<\beta_N$ given in Table \ref{table:summary}, and therefore  $\lambda^*(N) \le \lambda'_N$.
\item If $N=9$, then Lemma \ref{sing-lem} holds with $w:=w_{2.8}$, $\lambda'_9:=366<\beta_9:=368.5$, and therefore   $\lambda^*(9) \le 366$.
\end{enumerate}
The extremal solution is therefore singular for dimension $N\geq 9$.
\end{thm}

\begin{table}[ht]
\caption{Summary} 
\centering 
\begin{tabular}{c c c c} 
\hline\hline 
N&w & $\lambda'_{N}$ & $\beta_{N}$  \\ [0.5ex] 
\hline 
9 & $w_{2.8}$ &366 & 366.5  \\
10 & $w_3$ &450 & 487  \\
11& $w_3$ & 560 & 739  \\
12& $w_3$ & 680 & 1071 \\
13& $w_3$ & 802 & 1495 \\
14& $w_3$ & 940 & 2026 \\
15 & $w_3$& 1100 & 2678 \\
16 & $w_3$& 1260 & 3469  \\
$17\leq N \leq  30$& $w_3$& $H_N/2$ & $H_N/2$ \\
$N\geq  31$ & $w_2$& $27 {\bar \lambda}_N$ & $H_N/2$  \\
[1ex] 
\hline 

\end{tabular}
\label{table:summary} 
\end{table}
\begin{proof}
1)\, Assume first that  $N \ge 31$, then $27 \bar{ \lambda} \le
\frac{H_N}{2}$. We shall show that $w_2$ is a singular $H^2(B)-$weak
sub-solution of $(P)_{27 \bar \lambda}$ so that (\ref{cond2}) holds
with $\beta=\frac{H_N}{2}$.  Indeed,  write
\[
w_2:=1-|x|^{\frac{4}{3}}-2(|x|^\frac{4}{3}-|x|^2)=\bar{u}-\phi_0,
\]
 where $ \phi_0:=2(|x|^\frac{4}{3}-|x|^2)$, and note that  $w_2\in H_0^2(B)$, $ \frac{1}{1-w_2} \in L^3(B)$, $ 0 \le w_2 \le 1$ in $B$, and
\[
\Delta^2 w_2=\frac{3\bar \lambda}{r^{\frac{8}{3}}} \le \frac{ 27 \bar{\lambda}}{(1-w_2)^2} \qquad \hbox{in }B\setminus \{0\}.
\]
So $w_2$ is $H^2(B)-$weak sub-solution of $(P)_{27 \bar \lambda} $. Moreover, by $\phi_0\geq 0$ and (\ref{HR1}) we get that
\begin{eqnarray*}
H_N \int_B \frac{ \phi^2}{(1-w_2)^3} = H_N \int_B \frac{ \phi^2}{(|x|^\frac{4}{3}+\phi_0)^3} \le H_N \int_B \frac{\phi^2}{|x|^4} \le \int_B (\Delta \phi)^2
\end{eqnarray*}
for all $ \phi \in H_0^2(B)$. It follows from Lemma \ref{sing-lem} that $u^*$ is singular and that $ \lambda^* \leq 27 \bar{ \lambda} \le \frac{H_N}{2}$.\\
2)\,  Assume $ 17 \le N \le 30$ and consider the function
\[
w_3:=1- \frac{9}{5} r^\frac{4}{3} + \frac{4}{5} r^3.
\]
We show that $w_3$ is a semi-stable singular $H^2(B)-$weak sub-solution of $ (P)_\frac{H_N}{2}$. Indeed, we clearly have that  $ 0 \le w_3 \le 1$ in $B$, $ w_3 \in H_0^2(B)$ and $ \frac{1}{1-w_3} \in L^3(B)$.   To show the stability condition, we consider  $ \phi \in H^2_0(B)$ and write
\begin{eqnarray*}
H_N \int_B \frac{ \phi^2}{(1-w_3)^3} &=& 125 H_N \int_B \frac{ \phi^2}{ (9r^\frac{4}{3}-4r^3)^3} \le  125 H_N \sup_{0<r<1} \frac{1}{(9-4r^{\frac{5}{3}})^3} \int_B \frac{\phi^2}{r^4} \\
&=& H_N \int_B \frac{ \phi^2}{r^4} \le \int_B (\Delta \phi)^2
\end{eqnarray*}
by virtue of (\ref{HR1}). An easy computation shows that
\begin{eqnarray*}
\frac{H_N}{2(1-w_3)^2} - \Delta^2 w_3 &=& \frac{ 25 H_N}{2 ( 9r^\frac{4}{3}-4r^3)^2} - \frac{ 9 \bar{\lambda}}{5 r^\frac{8}{3}} - \frac{12}{5}\frac{N^2-1}{r}\\
&=& \frac{25 N^2 (N-4)^2 }{32 ( 9 r^\frac{4}{3}-4r^3)^2} - \frac{8 ( N-\frac{2}{3})(N-\frac{8}{3})}{ 5 r^\frac{8}{3}}  - \frac{12}{5}\frac{N^2-1}{r}.
\end{eqnarray*}
By using Maple one can verify that this final quantity is nonnegative on $(0,1)$ whenever  $ 17 \le N \le 30$, and hence $ w_3 $ is a $H^2(B)-$weak sub-solution of $ (P)_\frac{H_N}{2}$. It follows from Lemma \ref{sing-lem} that $u^*$ is singular and that $ \lambda^* \le \frac{H_N}{2}$.\\
3)\, Assume $10\leq N \leq 16$. We shall prove that again   $w:=w_{3}$ satisfies the assumptions of
Lemma \ref{sing-lem}. Indeed, using Maple,  we show that for each dimension $10\leq N\leq 16$, inequality (\ref{cond1}) holds with  $\lambda'_{N}$ given by Table \ref{table:summary}. Then, by using Maple again,  we show
that for each dimension $10\leq N\leq 16$, the following inequality holds
\begin{eqnarray*}
\frac{(N-2)^2(N-4)^2}{16}\frac{1}{(|x|^2-|x|^{\frac{N}{2}+1})(|x|^2-|x|^{\frac{N}{2}})}&+&\frac{(N-1)(N-4)^2}{4}
\frac{1}{|x|^2(|x|^2-|x|^{\frac{N}{2}})}\\ &\geq&
\frac{2\beta_{N}}{(1-w_{3})^3}.
\end{eqnarray*}
where $\beta_{N}$ is again given by Table \ref{table:summary}. The above inequality and the Hardy-Rellich inequality
(\ref{HR2}) guarantee that the stability condition (\ref{cond2})
holds with  $\beta:=\beta_{N}$. Since $\beta_{N}>\lambda'_{N}$, we deduce from Lemma \ref{sing-lem} that the
extremal solution is singular for $10\leq N\leq 16$.\\
4)\, Suppose now $N=9$ and consider $w:=w_{2.8}$. Using Maple on can see that
\[\Delta ^2 w \leq \frac{366}{(1-w)^2} \ \ \hbox{in} \ \ B\]
and
\[\frac{723}{(1-w)^3}\leq Q(r)\left(
P(r)+\frac{N-1}{r^2}\right) \ \ \hbox{for all} \ \ r\in(0,1), \]
where $P$ and $Q$ are given in  (\ref{HR3}). Since $723>2\times
366,$ by Lemma \ref{sing-lem} the extremal solution $u^{*}$ is
singular in dimension $N=9$.\end{proof}

\section{Appendix: Improved Hardy-Rellich Inequalities}
We now prove the improved Hardy-Rellich inequalities used in section 4. They rely on the results of Ghoussoub-Moradifam in \cite{GM} which provide necessary and sufficient conditions for such inequalities to hold. At the heart of this characterization is the following notion of a Bessel pair of functions.
\begin{dfn} Assume that $B$ is a ball of radius $R$ in $\IR^N$, $V,W \in C^{1}(0,1)$, and
$\int^{R}_{0}\frac{1}{r^{N-1}V(r)}dr=+\infty$. Say that the couple $(V, W)$  is a {\it Bessel pair on $(0, R)$} if the ordinary differential equation
\begin{equation*}
\hbox{ $({\rm B}_{V,W})$  \quad \quad \quad \quad \quad \quad \quad
\quad \quad \quad \quad
$y''(r)+(\frac{N-1}{r}+\frac{V_r(r)}{V(r)})y'(r)+\frac{W(r)}{V(r)}y(r)=0$
\quad \quad \quad \quad \quad \quad \quad \quad \quad \quad \quad}
\end{equation*}
has a positive solution on the interval $(0, R)$.
\end{dfn}
The space of radial functions in $C_{0}^\infty(B)$ will be denoted by $C^\infty_{0,r}(B)$.
The needed inequalities will follow from the following result.
\begin{thm}\label{GM1} {\bf (Ghoussoub-Moradifam \cite{GM}) } Let $V$ and $W$ be positive radial $C^1$-functions   on $B\backslash \{0\}$, where $B$ is a ball centered at zero with radius $R$ in $\IR^N$ ($N \geq 1$) such that  $\int^{R}_{0}\frac{1}{r^{N-1}V(r)}dr=+\infty$ and $\int^{R}_{0}r^{N-1}V(r)dr<+\infty$. The following statements are then equivalent:

\begin{enumerate}

\item $(V, W)$ is a Bessel pair on $(0, R)$.

\item $ \int_{B}V(|x|)|\nabla \phi |^{2}dx \geq \int_{B} W(|x|)\phi^2dx$ for all $\phi \in C^\infty_{0}(B)$.

\item If $\lim_{r \rightarrow 0}r^{\alpha}V(r)=0$ for some $\alpha< N-2$, then the above are
equivalent to
$$\int_{B}V(|x|)(\Delta \phi )^{2}dx \geq  \int_{B} W(|x|)|\nabla \phi|^{2}dx+(N-1)\int_{B}(\frac{V(|x|)}{|x|^2}-\frac{V_r(|x|)}{|x|})|\nabla \phi|^2dx$$
for all $\phi \in C^\infty_{0,r}(B)$.

\item If in addition, $W(r)-\frac{2V(r)}{r^2}+\frac{2V_r(r)}{r}-V_{rr}(r)\geq 0$ on $(0, R)$, then the above are equivalent to
$$\int_{B}V(|x|)(\Delta \phi)^{2}dx \geq  \int_{B} W(|x|)|\nabla \phi|^{2}dx+(N-1)\int_{B}(\frac{V(|x|)}{|x|^2}-\frac{V_r(|x|)}{|x|})|\nabla \phi|^2dx$$
for all  $\phi \in C^\infty_{0}(B)$.
\end{enumerate}

\end{thm}
We shall now deduce the following corollary.
\begin{cor}\label{GM2}
Let $N\geq 5$ and $B$ be the unit ball in $\IR^N$. Then the
following improved Hardy-Rellich inequality holds for all $\phi \in
C^\infty_{0}(B)$:
\begin{eqnarray}\label{HR11}
\int_{B}(\Delta \phi)^2 &\geq&
\frac{(N-2)^2(N-4)^2}{16}\int_{B}\frac{\phi^2}{(|x|^2-|x|^{\frac{N}{2}+1})(|x|^2-|x|^{\frac{N}{2}})}\\
&&+\frac{(N-1)(N-4)^2}{4} \int_{B}\frac{\phi^2}{|x|^2(|x|^2-|x|^{\frac{N}{2}})}\nonumber
\end{eqnarray}
\end{cor}
\begin{proof} Let $0<\alpha<1$ and define $y(r):=r^{-\frac{N}{2}+1}-\alpha$. Since
\begin{equation*}
-\frac{y''+\frac{(N-1)}{r}y'}{y}=\frac{(N-2)^2}{4}\frac{1}{r^2-\alpha
r^{\frac{N}{2}+1}},
\end{equation*}
the couple $\left(1,\frac{(N-2)^2}{4}\frac{1}{r^2-\alpha
r^{\frac{N}{2}+1}}\right)$ is a Bessel pair on $(0,1)$. By Theorem
\ref{GM1}(4) the following inequality then holds:
\begin{equation}\label{AHR}
\int_{B}(\Delta \phi)^2 dx\geq \frac{(N-2)^2}{4}
\int_{B}\frac{|\nabla \phi|^2}{|x|^2-\alpha
|x|^{\frac{N}{2}+1}}+(N-1)\int_{B}\frac{|\nabla \phi|^2}{|x|^2}
\end{equation}
for all $\phi \in C^\infty_0(B)$. Set $V(r):=\frac{1}{r^2-\alpha
r^{\frac{N}{2}+1}}$ and note that
$$\frac{V_{r}}{V}=-\frac{2}{r}+\frac{\alpha(N-2)}{2} \frac{r^{\frac{N}{2}-2}}{1-\alpha r^{\frac{N}{2}-1}}\geq-\frac{2}{r}.$$
The function $y(r)=r^{-\frac{N}{2}+2}-1$ is decreasing and is then a
positive super-solution on $(0,1)$ for the ODE
\begin{equation*}
y''+(\frac{N-1}{r}+\frac{V_{r}}{V})y'(r)+\frac{W_{1}(r)}{V(r)}y=0,
\end{equation*}
where
\[W_{1}(r)=\frac{(N-4)^2}{4(r^2-r^{\frac{N}{2}})(r^2-\alpha r^{\frac{N}{2}+1})}.\]
Hence, by Theorem \ref{GM1}(2) we deduce
$$\int_{B}\frac{|\nabla
\phi|^2}{|x|^2-\alpha|x|^{\frac{N}{2}+1}}\geq(\frac{N-4}{2})^2\int_{B}\frac{\phi^2}{(|x|^2-\alpha
|x|^{\frac{N}{2}+1})(|x|^2-|x|^{\frac{N}{2}})}$$ for all $\phi \in
C_0^\infty(B)$. Similarly, for $V(r)=\frac{1}{r^2}$ we have that
$$\int_{B}\frac{|\nabla \phi|^2}{|x|^2}\geq
(\frac{N-4}{2})^2\int_{B}\frac{\phi^2}{|x|^2(|x|^2-|x|^{\frac{N}{2}})}$$
for all $\phi \in C_0^\infty(B)$. Combining the above two
inequalities with (\ref{AHR}) and letting $\alpha \rightarrow 1$ we
get inequality (\ref{HR11}).
\end{proof}

\begin{cor}\label{GM3}
Let $N=9$ and $B$ be the unit ball in $\IR^N$. Define
$\varphi(r):=r^{-\frac{N}{2}+1}+r-1.9$ and
$\psi(r):=r^{-\frac{N}{2}+2}+20r^{-1.69}+10r^{-1}+10r+7r^2-48$. Then
the following improved Hardy-Rellich inequality holds for all $\phi
\in C^{\infty}_{0}(B)$:
\begin{equation}\label{HR13}
\int_{B}(\Delta \phi)^2\geq \int_{B}Q(|x|)\left(
P(|x|)+\frac{N-1}{|x|^2}\right)\phi^2,
\end{equation}
where
\[P(r):=-\frac{\varphi''(r)+\frac{N-1}{r}\varphi'(r)}{\varphi(r)} \ \ \hbox{and} \ \ Q(r):=-\frac{\psi''(r)+\frac{N-3}{r}\psi'(r)}{\psi(r)}. \]
\end{cor}
\begin{proof} By definition  $(1,P(r))$ is a Bessel pair on
$(0,1)$. One can easily see that $P(r)\geq \frac{2}{r^2}$. Hence, by
Theorem \ref{GM1}(4) the following inequality holds:
\begin{equation}\label{AHRbis}
\int_{B}(\Delta \phi)^2 dx\geq  \int_{B}P(|x|)|\nabla
\phi|^2+(N-1)\int_{B}\frac{|\nabla \phi|^2}{|x|^2}
\end{equation}
for all $\phi \in C^{\infty}_{0}(B)$. Using Maple it is easy to see that
$$\frac{P_{r}}{P}\geq-\frac{2}{r}\ \ \hbox{in} \ \ (0,1),$$
and therefore $\psi(r)$ is a positive super-solution for the ODE
\begin{equation*}
y''+(\frac{N-1}{r}+\frac{P_{r}(r)}{P(r)})y'(r)+\frac{P(r)Q(r)}{P(r)}y=0,
\end{equation*}
on $(0,1)$.  Hence, by Theorem \ref{GM1}(2) we have for all $\phi
\in C_0^\infty(B)$
$$\int_{B}P(|x|)|\nabla \phi|^2 \geq\int_{B}P(|x|)Q(|x|)\phi^2,$$
and similarly
$$\int_{B}\frac{|\nabla \phi|^2}{|x|^2}\geq
\int_{B}\frac{Q(|x|)}{|x|^2}\phi^2,$$
since $\psi(r)$ is a positive solution for the ODE
\begin{equation*}
y''+\frac{N-3}{r}y'(r)+Q(r)y=0.
\end{equation*}
Combining the above two inequalities with (\ref{AHRbis}) we get
(\ref{HR13}).
\end{proof}

\end{document}